# A new signal processing tool developed with the help of the Clifford algebra

Amitabha Chanda

Amitabha Chanda, Guest Professor, UCSTA, University of Calcutta,
92 APC Road, Kolkata 700009
Communication address: Amitabha Chanda, P36/1, Manasbag, Kolkata 700056, India.
E-mail: amitabha39@yahoo.co.in
amitabha39@gmail.com

Phones: 0091 98308 55481 (M), (0091) (33) 2553 7298 (R)

## ABSTRACT

A positive integer is expressed as a sum of squares of positive integers in a unique way applying a special technique. The expression, thus obtained is resolved into two factors using the concept of Clifford algebra. This technique is applied on a set of positive integers. Using these factors, a new mean of the set of integers, different from the arithmetic mean, is developed. A detailed discussion involving Cl(0,3) as an example is presented. The new mean bears the imprint of degree of randomness of the data set. More random is the data set more away is the new mean from the arithmetic mean. In addition to the new mean a new standard deviation is also developed. It is shown that the new standard deviation is more sensitive to randomness than standard deviation.

# 1. Introduction

Statistical analysis is an essential tool of almost all branches of knowledge. It is adequately dealt in the books of Yale and Kendal [1], Weatherburn [2] and many others. The fundamentals of statistical analysis are distribution functions, calculation of moments. The calculation of moments needs elementary arithmetical operation like arithmetic averaging. In the present paper we develop an averaging technique, based on **Clifford algebra**. It is altogether different from arithmetic averaging. .

In addition to Clifford algebra we use technique, developed for the present purpose. In this technique a positive integer is expressed as a sum of perfect squares in a unique way. A name is given to this novel technique: **SSGS**. The combination of **SSGS** and **Clifford algebra** is used in calculating new types of mean and standard deviation. .

The application of the suggested technique to some positive decimals is also discussed.

# 2. A representation of positive integers

Besides decimal and binary number systems there exist many ways to represent real numbers. For example the expansion of numbers in the form of a continued fraction may be cited. The conventional representation of real numbers in base **b, b ≥ 2**, has been generalized by A. R´enyi [3]. He considered writing non-negative real number **x** with arbitrary real base **β > 1**.

Lagrange showed in his famous Four Square Theorem that a positive integer can be expressed as a sum of only four perfect squares.. But for a particular positive integer there may be more than one representation. An example may explain this.

$$191 = 10^2 + 9^2 + 3^2 + 1^2 = 11^2 + 6^2 + 5^2 + 3^2$$

Chanda [4] introduces a new system of representing positive integers. It is observed that with this system the representation is unique. We shall first discuss the method and then in the subsequent section apply it to unravel some basic statistical properties of some positive integers and positive decimals. For that purpose we shall take help of Clifford Algebra.

**Comment**

*For each nonnegative integer **x** there exists a nonnegative integer $x_1$, such that*

$$x_1^2 \leq x < (x_1 + 1)^2 \qquad (1)$$

**A few definitions**

(a) **Significant square (ss)** in a given nonnegative integer **x** is the square of nonnegative integer which is just less or equal to **x.** As in (1) $x_1^2$ is the significant square, extracted from **x.** (1) may be rewritten as

$$x_1^2 \leq x \leq x_1^2 + 2x_1 \qquad (2)$$

(b) **Residue $R_1$** is given as

$$R_1 = x - x_1^2 \qquad (3)$$

$R_1$ is again a nonnegative integer.

(c) First **generation square (gs)** $x_2^2$ is given as

$$x_2^2 \leq R_1 \leq x_2^2 + 2x_2 \qquad (4)$$

In fact $x_2$ is the greatest positive integer whose square is less than or equal to $R_1$.

*The process goes on like this till the residue is zero. In a trivial case if the original x is a perfect square then there is no non zero residue. Following this process a nonnegative integer **x** can be expressed as*

$$x = x_1^2 + x_2^2 + \cdots\cdots + x_n^2 \qquad (5)$$

*Here $x_1^2$ is the significant square and $x_2^2, ... x_n^2$ are generation squares, as defined earlier in.*

**Some examples**

i) $91 = 9^2 + 3^2 + 1^2$

Here $9^2$ is significant square (ss) and $3^2$ and $1^2$ are generation squares (gs).

(ii) $192 = 13^2 + 4^2 + 2^2 + 1^2 + 1^2 + 1^2$

Here $13^2$ is significant square (ss) and $4^2$, $2^2$, $1^2$, $1^2$ and $1^2$ are generation squares (gs).

(iii) $999998 = 999^2 + 44^2 + 7^2 + 3^2 + 1^2 + 1^2 + 1^2$

## 3. A brief review of the Clifford Algebra

At the outset we introduce in short the algebra, developed by Clifford [5]. The properties of Clifford numbers vis-à-vis algebra have been elaborately studied by Keller [6], Lounesto [7], Ablamowicz [8], Hestenes & Sobczyk [9], Sobczyk [10] and Oziewicz [11] and many others. For different types of Clifford numbers there are different types of bases. Sobczyk [10] has discussed different types of bases for $R^2$ and $R^3$ When the number of bases is **n** the Clifford number will have $2^n$ elements. We would confine our application to Clifford numbers **Cl(0,3)**. A sample pair of **Cl(0,3)** number may be given as below.

$$\begin{aligned} p &= \varepsilon_0 p_0 + \varepsilon_1 p_1 + \varepsilon_2 p_2 + \varepsilon_3 p_3 + \varepsilon_4 p_4 + \varepsilon_5 p_5 + \varepsilon_6 p_6 + \varepsilon_7 p_7' \\ p^c &= \varepsilon_0 p_0 - \varepsilon_1 p_1 - \varepsilon_2 p_2 - \varepsilon_3 p_3 - \varepsilon_4 p_4 - \varepsilon_5 p_5 - \varepsilon_6 p_6 - \varepsilon_7 p_7' \end{aligned} \qquad (7)$$

The first one is the intended Clifford number and the second, in absence of a suitable term we call it as the conjugate of the first one. . The bases in (7) have the following properties.

$$\varepsilon_4 = \varepsilon_1 \wedge \varepsilon_2 \quad \varepsilon_5 = \varepsilon_2 \wedge \varepsilon_3 \quad \varepsilon_6 = \varepsilon_3 \wedge \varepsilon_1 \quad \varepsilon_7 = \varepsilon_1 \wedge \varepsilon_2 \wedge \varepsilon_3 \tag{8}$$

$$\varepsilon_0^2 = +1 \quad \varepsilon_1^2 = \varepsilon_2^2 = \varepsilon_3^2 = \varepsilon_4^2 = \varepsilon_5^2 = \varepsilon_6^2 = -1 \quad \varepsilon_7^2 = +1 \tag{9}$$

$p_k$'s [k=0…7] in (7) may be any type of numbers. But now we are only considering positive integers with one exception, $p_7'$ which is imaginary and is given as

$$p_7' = ip_7 \qquad i = \sqrt{-1}$$

Accordingly (7) may be rewritten as

$$\begin{aligned} p &= \varepsilon_0 p_0 + \varepsilon_1 p_1 + \varepsilon_2 p_2 + \varepsilon_3 p_3 + \varepsilon_4 p_4 + \varepsilon_5 p_5 + \varepsilon_6 p_6 + \varepsilon_7 (ip_7) \\ p^c &= \varepsilon_0 p_0 - \varepsilon_1 p_1 - \varepsilon_2 p_2 - \varepsilon_3 p_3 - \varepsilon_4 p_4 - \varepsilon_5 p_5 - \varepsilon_6 p_6 - \varepsilon_7 (ip_7) \end{aligned} \tag{10}$$

Because of the presence of an imaginary number **($ip_7$), ($p_7$ being a nonnegative integer)**, as the coefficient of the trivector, we would like to call it a **Cl(0,3)** number only  For developing these Clifford numbers the help of the reputed works of Ablamowitz [12] and Oziewicz [13] has been taken.

## 4. Application of the Clifford Algebra to the SSGS system of nonnegative integers

We have already observed in the foregoing paragraphs that following the **SSGS** system a very large nonnegative integer can be expressed as the sum of a few perfect squares and that also in a unique way. In view of the above we would have the liberty to carry out further study on sufficiently large positive integers even restricting their expressions to eight squares. It has already been exhibited by the above examples.

We shall now prove a theorem involving nonnegative integers and its Clifford representation.

**Theorem**

*The arithmetic average of a set of nonnegative integers, expressed as per the SSGS scheme as a sum of squares of positive integers is greater or equal to the product of two mutually conjugate Clifford numbers, the coefficients of which are the arithmetic averages of the positive integers, constituting the nonnegative integers..*

**Proof :**

Let us now arbitrarily take a set of **n** number of nonpositive integers each of which fits the restriction of eight squares. Any of these integers $a_i$ as per the SSGS scheme may be expressed as below.

$$a_i = p_{i0}^2 + p_{i1}^2 + p_{i2}^2 + p_{i3}^2 + p_{i4}^2 + p_{i5}^2 + p_{i6}^2 + p_{i7}^2 \tag{11}$$

Following the rules of Clifford algebra $\mathbf{a}_i$ may now be resolved into two mutually conjugate factors.

$$\mathbf{a}_i = \mathbf{p}_i \mathbf{p}_i^c \tag{12}$$

Here $\mathbf{p}_i$ is a **Cl(0,3)** Clifford number $\mathbf{p}_i^c$ is its conjugate similar to the expression (10).

We rewrite (12) as per the SSGS scheme explained earlier.

$$\begin{aligned}\mathbf{a}_i &= \mathbf{p}_{i0}^2 + \mathbf{p}_{i1}^2 + \mathbf{p}_{i2}^2 + \mathbf{p}_{i3}^2 + \mathbf{p}_{i4}^2 + \mathbf{p}_{i5}^2 + \mathbf{p}_{i6}^2 + \mathbf{p}_{i7}^2 \\ &= \{\varepsilon_0 \mathbf{p}_{i0} + \varepsilon_1 \mathbf{p}_{i1} + \varepsilon_2 \mathbf{p}_{i2} + \varepsilon_3 \mathbf{p}_{i3} + \varepsilon_4 \mathbf{p}_{i4} + \varepsilon_5 \mathbf{p}_{i5} + \varepsilon_6 \mathbf{p}_{i6} + \varepsilon_7 (i\mathbf{p}_{i7})\} \\ &\quad \times \{\varepsilon_0 \mathbf{p}_{i0} - \varepsilon_1 \mathbf{p}_{i1} - \varepsilon_2 \mathbf{p}_{i2} - \varepsilon_3 \mathbf{p}_{i3} - \varepsilon_4 \mathbf{p}_{i4} - \varepsilon_5 \mathbf{p}_{i5} - \varepsilon_6 \mathbf{p}_{i6} - \varepsilon_7 (i\mathbf{p}_{i7})\}\end{aligned} \tag{13}$$

As per the scheme the $\mathbf{p}_{ij}$'s in (13) are positive integers. It may be noted that the last few squares may be $\mathbf{0}^2$ also. Even some $\mathbf{a}_i$ may be a perfect square with only one term in rhs of (13).

Now summing up all the $\mathbf{a}_i$'s we get

$$\sum_{i=1}^n \mathbf{a}_i = \sum_{i=1}^n \mathbf{p}_i \mathbf{p}_i^c = \sum_{i=1}^n \mathbf{p}_i \sum_{j=1}^n \mathbf{p}_i^c - \sum_{i=1}^n \mathbf{p}_i \sum_{\substack{j=1 \\ j \neq i}}^n \mathbf{p}_i^c \tag{14}$$

The second term in lhs of (13) may be written as below.

$$\sum_{i=1}^n \mathbf{p}_i \sum_{\substack{j=1 \\ j \neq i}}^n \mathbf{p}_i^c = \frac{n-1}{n} \sum_{i=j}^n \mathbf{p}_i \sum_{j=1}^n \mathbf{p}_i^c - n\lambda \tag{15}$$

In (15) $\mathbf{n}\lambda$ is a number to be explained in the subsequent paragraphs. From (14) and (15) we get

$$\begin{aligned}\sum_{i=1}^n \mathbf{a}_i &= \sum_{i=1}^n \mathbf{p}_i \sum_{j=1}^n \mathbf{p}_i^c - \frac{n-1}{n} \sum_{i=j}^n \mathbf{p}_i \sum_{j=1}^n \mathbf{p}_i^c + n\lambda \\ &= \frac{1}{n} \sum_{i=1}^n \mathbf{p}_i \sum_{j=1}^n \mathbf{p}_i^c + n\lambda\end{aligned} \tag{16}$$

{vide Appendix A for detailed calculation of (13), (14) & (15)].

Now the arithmetic mean AM of $\mathbf{a}_i$'s is given as

$$\begin{aligned}\frac{1}{n}\sum_{i=1}^n \mathbf{a}_i &= \frac{1}{n^2}\sum_{i=1}^n \mathbf{p}_i \sum_{j=1}^n \mathbf{p}_i^c + \lambda \\ &= \left(\frac{1}{n}\sum_{i=1}^n \mathbf{p}_i\right)\left(\frac{1}{n}\sum_{j=1}^n \mathbf{p}_i^c\right) + \lambda\end{aligned} \tag{17}$$

Using (13), (17) is rewritten as below.

$$\frac{1}{n}\sum_{i=1}^{n}\left(\begin{array}{l}p_{i0}^{2}+p_{i1}^{2}+p_{i2}^{2}+p_{i3}^{2}+p_{i4}^{2}+p_{i5}^{2}\\+p_{i6}^{2}+p_{i7}^{2}\end{array}\right)$$

$$-\left(\frac{1}{n}\sum_{i=1}^{n}\{\varepsilon_{0}p_{i0}+\varepsilon_{1}p_{i1}+\varepsilon_{2}p_{i2}+\varepsilon_{3}p_{i3}+\varepsilon_{4}p_{i4}+\varepsilon_{5}p_{i5}+\varepsilon_{6}p_{i6}+\varepsilon_{7}(ip_{i7})\}\right)$$

$$\left(\frac{1}{n}\sum_{j=1}^{n}\{\varepsilon_{0}p_{i0}-\varepsilon_{1}p_{i1}-\varepsilon_{2}p_{i2}-\varepsilon_{3}p_{i3}-\varepsilon_{4}p_{i4}-\varepsilon_{5}p_{i5}-\varepsilon_{6}p_{i6}-\varepsilon_{7}(ip_{i7})\}\right)=\lambda$$

Or

$$\left(\overline{p_0^2}+\overline{p_1^2}+\overline{p_2^2}+\overline{p_3^2}+\overline{p_4^2}+\overline{p_5^2}+\overline{p_6^2}+\overline{p_7^2}\right)-\left(\overline{p_0}^2+\overline{p_1}^2+\overline{p_2}^2+\overline{p_3}^2+\overline{p_4}^2+\overline{p_5}^2+\overline{p_6}^2+\overline{p_7}^2\right)=\lambda$$

Or

$$\left(\overline{p_0^2}-\overline{p_0}^2\right)+\left(\overline{p_1^2}-\overline{p_1}^2\right)+\left(\overline{p_2^2}-\overline{p_2}^2\right)+\left(\overline{p_3^2}-\overline{p_3}^2\right)+\left(\overline{p_4^2}-\overline{p_4}^2\right)$$
$$+\left(\overline{p_5^2}-\overline{p_5}^2\right)+\left(\overline{p_6^2}-\overline{p_6}^2\right)+\left(\overline{p_7^2}-\overline{p_7}^2\right)=\lambda \quad (18)$$

In (18) each of the eight items within first bracket in l.h.s. is equal to the variance of the specific coefficient and hence a nonnegative quantity. So the item in r.h.s. is also nonnegative. Accordingly we may write the following inequality.

$$\frac{1}{n}\sum_{i=1}^{n}a_i \geq \left(\frac{1}{n}\sum_{i=1}^{n}p_i\right)\left(\frac{1}{n}\sum_{j=1}^{n}p_i^{\,c}\right) \quad (19)$$

The rhs of (19) is the product of the mean of Clifford numbers and the mean of their conjugates, and hence a real number. We shall give it a name $\mu_{1c}$, expressed as below.

$$\mu_{1c}=\left(\frac{1}{n}\sum_{i=1}^{n}p_i\right)\left(\frac{1}{n}\sum_{j=1}^{n}p_i^{\,c}\right) \quad (20)$$

Finally we get the relation from (17), (18) and (19) as below.

$$\overline{a}=\mu_{1c}+\lambda$$
$$\overline{a}\geq\mu_{1c} \quad (21)$$

QED

## 5. Extension of the scheme to finite positive decimals

A positive decimal number with n decimal places may be expressed as below.

$$a = a_m 10^m + a_{m-1} 10^{m-1} + \cdots\cdots + a_1 10^1$$
$$+ a_0 10^0 + \cdots\cdots + a_{-n+1} 10^{-n+1} + a_{-n} 10^{-n} \qquad (27)$$

Here **n = 2k** or **2k-1** i.e even or odd. To proceed further let us make **n** an even number. It can be done in two ways either by deleting the number in the last decimal place following laws of approximation or by adding a zero at the last to make **n=2k**. We shall choose the second one since by this we neither miss any information related to the number, as we do if we follow the first one; nor run the risk of having different **SSGS** expressions for different number of pairs of zeroes. We shall explain it later; but for the time being let us proceed by taking **n=2k**, if necessary by addition of one zero at the end. Then we get the expression of a as below.

$$a = a_m 10^m + a_{m-1} 10^{m-1} + \cdots\cdots + a_1 10^1 + a_0 10^0$$
$$+ \cdots\cdots + a_{-n+1} 10^{-2k+1} + a_{-n} 10^{-2k} \qquad (28)$$

The expression (27) may be expressed as below.

$$a = \frac{(a_m 10^{m+2k} + a_{m-1} 10^{m+2k-1} + \cdots + a_0 10^{2k} + \cdots + a_{-n+1} 10 + a_{-n} 10^0)}{10^{2k}} \qquad (29)$$

The numerator of (29) is a positive integer. Following **SSGS** system it is possible to express the numerator as a sum of squares of eight squares. In case there are more than eight squares **Cl(0,n)** for suitable **n** will be selected for further processing. However, we are restricting our discussion to **Cl(0,3)**, as in the case of positive integer. Then we get from (29)

$$a = \frac{q_0^2 + q_1^2 + q_2^2 + q_3^2 + q_4^2 + q_5^2 + q_6^2 + q_7^2}{10^{2k}}$$
$$= \left(\frac{q_0}{10^k}\right)^2 + \left(\frac{q_1}{10^k}\right)^2 + \left(\frac{q_2}{10^k}\right)^2 + \left(\frac{q_3}{10^k}\right)^2 + \left(\frac{q_4}{10^k}\right)^2 + \left(\frac{q_5}{10^k}\right)^2 + \left(\frac{q_6}{10^k}\right)^2 + \left(\frac{q_7}{10^k}\right)^2 \qquad (30)$$

In (30) we have expressed the number **a** as a sum of squares (not integers) as per **SSGS** system. From this point we may proceed as we do for positive integers.

NB. Let us take a case of odd **n**.

$$12.3 = 12.30 = \frac{1230}{100} = \frac{35^2 + 2^2 + 1^2}{100} = (3.5)^2 + (0.2)^2 + (0.1)^2$$

$$12.3 = 12.3000 = \frac{123000}{10000} = \frac{(350)^2 + (22)^2 + (4)^2}{10000} = (3.5)^2 + (0.22)^2 + (0.4)^2$$

We shall restrict our operation to the first one, i.e. add only one zero at the end to make it same as the case of even **n.** We would avoid the second expression since we get a different result due to a few zeroes that do not influence the input.

## 6. Properties of the New Mean of the set of positive integers and also of positive decimals as discussed

      i. New Mean is always less than arithmetic mean until and unless the data set is completely nonrandom.
      ii. More random is the data set more is the difference between New Mean and arithmetic mean.
      iii. Unlike arithmetic mean New Mean is sensitive to the distribution pattern of the set.

In the subsequent paragraph we would give an example in support of our finding.

**Example**

**Set 1: 101, 118, 99, 131, 140, 141, 109, 121, 122, 110. Total= 1192 AM= 119.2**
**Set 2: 112, 107, 103, 135, 131, 130, 123, 109, 130, 112. Total = 1192 AM= 119.2**

It may be noted that in both the sets total and AM's are same. The results of Set 1 and Set 2 as per the above process are given as below.

$$\overline{q_0^2} - \overline{q_0}^2 = 0.44 \quad \overline{q_1^2} - \overline{q_1}^2 = 2.01 \quad \overline{q_2^2} - \overline{q_2}^2 = 0.41$$

$$\overline{q_3^2} - \overline{q_3}^2 = 0.21 \quad \overline{q_4^2} - \overline{q_4}^2 = 0.09 \quad \overline{q_5^2} - \overline{q_5}^2 = \overline{q_6^2} - \overline{q_6}^2 = \overline{q_7^2} - \overline{q_7}^2 = 0$$

**Total variance = 1.72   New Mean = 117.48 AM= 117.48+1.72= 119.2**

**Calculation of Set 1.**

$$101 = 10^2 + 1^2 \quad 118 = 10^2 + 4^2 + 1^2 + 1^2 \quad 99 = 9^2 + 4^2 + 1^2 + 1^2$$

$$131 = 11^2 + 3^2 + 1^2 \quad 140 = 11^2 + 4^2 + 1^2 + 1^2 + 1^2 \quad 141 = 11^2 + 4^2 + 2^2$$

$$109 = 10^2 + 3^2 \quad 121 = 11^2 \quad 122 = 11^2 + 1^2 \quad 110 = 10^2 + 3^2 + 1^2$$

New Mean is then calculated as below. We give the detailed calculation for Set 1.

$$\left\{\frac{1}{10}(10+10+9+11+11+11+10+11+11+10)\right\}^2 + \left\{\frac{1}{10}(1+4+4+3+4+4+3+1+3)\right\}^2$$

$$+ \left\{\frac{1}{10}(1+1+1+1+2+1)\right\}^2 + \left\{\frac{1}{10}(1+1+1)\right\}^2 + \left\{\frac{1}{10}(1)\right\}^2 = 116.04$$

**Set 2 :**  By similar calculation we get **New Mean = 117.48**

## 7. Application of the New Mean in signal processing

In the subsequent example we have made an attempt to explain the usefulness of the New Mean. Signal processing data, realized from a continuous signal at a fixed interval, may be expressed as a time series. We scale the data in such a manner that all the data become positive integers with the minimum data being **1** Then we proceed to find the MA taking a particular length of string. For each set the New Mean is calculated. We call these averages as the New Moving Average (NMA). The salient difference between the MA and the NMA is found out from two sets of data selected from the time series, where the sums of the respective data are same and so also the AM's. It is shown by calculating the NM's , the ranges and the SD,s;  that the gap between the AM and the NM increases as the randomness increases. The final results are shown in the fig.1 and fig.2 given below. In fig.1 an arbitrarily built signal with all data in positive integers is presented. In fig. 2 standard MA with 5 data length string is plotted against NMA with the same data string.

| Position of average | 36 | 49 |
|---|---|---|
| Length of string | 5 | 5 |
| Total of data | 717 | 717 |
| AM | 143.4 | 143.4 |
| New Mean | 141.88 | 140.76 |
| Range | 43 | 88 |
| Standard deviation | 26.11 | 36.64 |
| Data | 128,128,145,145,171 | 156,107,127,195,182 |

**Table 1. A comparative study between two strings having equal data value total; but different distributions.**

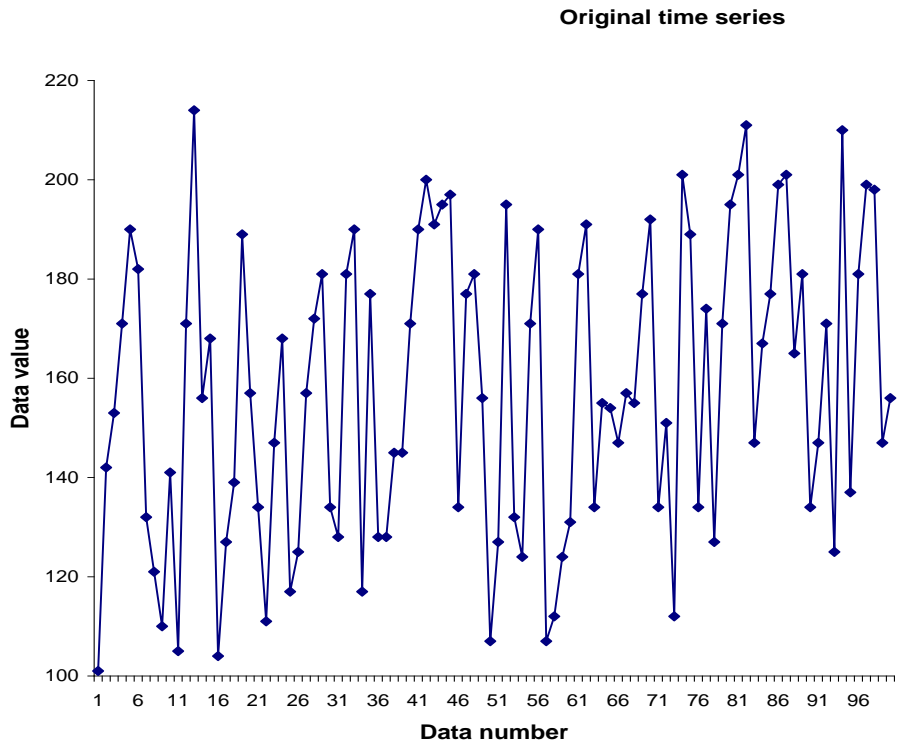

**Fig.1. An arbitrarily chosen time series.**

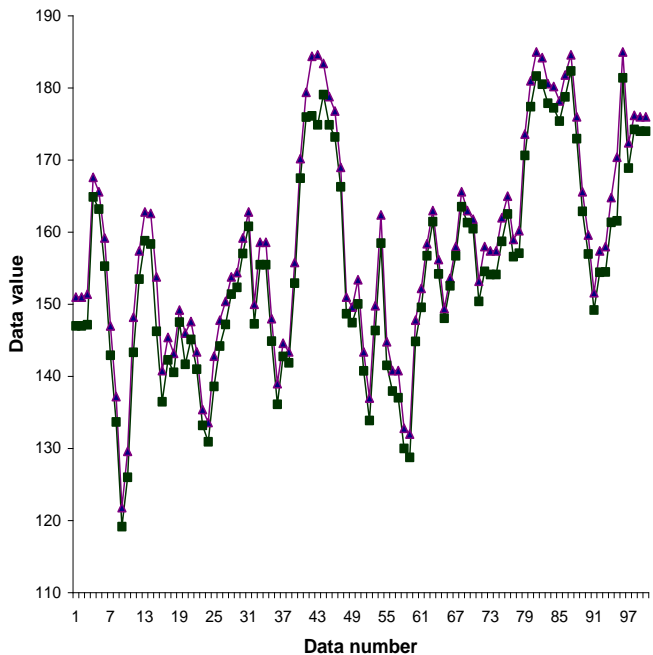

**Fig.2. Moving average (MA) vs. new moving average (NMA). Triangles show MA and squares NMA. Both have been built from data in fig.1 with string length 5**

## 9.Conclusion

We have developed and discussed a new mean. We have also tried to show its salient features in the background of the properties of arithmetic mean. It is observed that the new mean., developed here is sensitive to the randomness of the data set. In fact if there is no randomness both the arithmetic mean and the new mean are same.

**Acknowledgement**

I acknowledge the valuable advice I received from Professor Zbigniew Oziewicz, the eminent mathematician in the field of Clifford algebra.

**Appendix A:** Details of equation **(15), (16) & (17)**

$$\sum_{i=1}^{n} a_i$$

$$= \sum_{i=1}^{n} p_i p_i^{\;c} = \left(p_1 p_1^{\;c} + p_2 p_2^{\;c} + \cdots\cdots + p_n p_n^{\;c}\right)$$

$$= \left(p_1 + p_2 + \cdots\cdots + p_n\right)\left(p_1^{\;c} + p_2^{\;c} + \cdots\cdots + p_n^{\;c}\right)$$

$$-\left\{p_1\left(p_2^{\;c} + p_3^{\;c} + \cdots\cdots + p_n^{\;c}\right) + p_2\left(p_1^{\;c} + p_3^{\;c} + \cdots\cdots + p_n^{\;c}\right)\right.$$

$$\left.+ \cdots\cdots + p_n\left(p_1^{\;c} + p_2^{\;c} + \cdots\cdots + p_{n-1}^{\;c}\right)\right\}$$

$$= \sum_{i=1}^{n} p_i \sum_{j=1}^{n} p_j^{\;c} - \sum_{i=1}^{n} p_i \sum_{\substack{j=1 \\ j \neq i}}^{n} p_j^{\;c}$$

$$= \sum_{i=1}^{n} p_i \sum_{j=1}^{n} p_j^{\;c} - \frac{n(n-1)}{n^2} \sum_{i=1}^{n} p_i \sum_{j=1}^{n} p_j^{\;c} + n\lambda$$